\documentclass{article}
\usepackage{amsmath, amsthm, amssymb, geometry, color}
\theoremstyle{definition}
 \newtheorem{thm}{Theorem}
 
 \newtheorem{rem}{Remark}

\begin{document}

\title{Spectra of the Rarita-Schwinger operator on some symmetric spaces}
\author{Yasushi Homma, Takuma Tomihisa}
\date{}
\maketitle

\renewcommand{\thefootnote}{}
\footnotetext{{\em keyword} : Dirac operator, Rarita-Schwinger operator, Casimir operator on symmetric spaces}
\renewcommand{\thefootnote}{\arabic{footnote}}

\begin{abstract}
We give a method to calculate spectra of the square of the Rarita-Schwinger operator on compact symmetric spaces. According to Weitzenb\"{o}ck formulas, the operator can be written by the Laplace operator, which is the Casimir operator on compact symmetric spaces. Then we can obtain the spectra by using the Freudenthal's formula and branching rules. As examples, we calculate the spectra on the sphere, the complex projective space, and the quaternionic projective space.
\end{abstract}

\section{Introduction}

The Rarita-Schwinger operator is a generalization of the classical Dirac operator. The Dirac operator is a famous first order elliptic differential operator on spin \(\frac{1}{2}\) fields and has been studied for a long time not only in physics but also in analysis and differential geometry. The Rarita-Schwinger operator is ``spin \(\frac{3}{2}\) version'' of  the Dirac operator. There has been much research about the Rarita-Schwinger operator on the Euclidean spaces from the view point of Clifford analysis. For example, polynomial solutions, fundamental solution and Clifford-Cauchy kernel were studied in \cite{BSSL1}, \cite{BSSL2}. Furthermore, the generalization to the higher spin operators is in progress, \cite{ERJ}, etc. On non-Euclidean spaces, the harmonic analysis of the Dirac operator has been studied a lot. The eigenvalues of the Dirac operator was calculated explicitly on the sphere (cf. \cite{Bar}), the odd dimensional complex projective space \cite{CG}, \cite{SS}, the quaternionic projective space \cite{M} and more many examples, for more details see \cite{Gi}. On the other hand, there are a few results about the harmonic analysis of the Rarita-Schwinger operator on non-Euclidean spaces. In fact, the eigenvalues of the Rarita-Schwinger operator was calculated only in two cases.  One is the sphere case in \cite{Br}, \cite{BS}, the other is \(S^1\times S^n\) case in \cite{Do}. One of the reason why we could have done only them is that it is difficult to apply the method of calculation on the sphere in \cite{Br}, \cite{BS} (originally from \cite{BOO}) to other symmetric spaces. 

 Recently, the Rarita-Schwinger operators on Einstein manifolds were researched in \cite{HS}. In the paper, they show that the Rarita-Schwinger operator is the linear combination of the scalar curvature and the standard Laplacian on a Einstein manifold by using the Weitzenb\"{o}ck formulas. Here, the standard Laplacian is one of important operators naturally defined on Riemannian and spin manifolds \cite{H}, \cite{SW}. Especially, it coincides with the Casimir operator on compact symmetric spaces.

In this paper, we give a way to calculate the eigenvalues of the Rarita-Schwinger operator on irreducible compact symmetric spaces, which are famous examples of Einstein manifolds. Because the twistor operator \(P\) is overdetermined elliptic, the space of the spin \(\frac{3}{2}\) fields is decomposed into the direct sum of the kernel of \(P^\ast\) and the image of \(P\). According to the Weitzenb\"{o}ck formula, we write the square of the Rarita-Schwinger operator in terms of the Casimir operator and the scalar curvature on each components. Therefore, we can arrive at the eigenvalues of the Rarita-Schwinger operator by using the Freudenthal's formula and branching rules. Our method is much simpler than before \cite{Br}, \cite{Do}. Moreover, we apply our method to some concrete examples. We actually calculate the eigenvalues for the sphere (theorem \ref{sphere}), the complex projective space (theorem \ref{cplx}) and the quaternionic projective space (theorem \ref{qtn}). These are all irreducible compact rank-1 symmetric spaces except for the octonionic projective plane \(\mathrm{F}_4/\mathrm{Spin}(9)\).\\

Acknowledgement: this research was partially supported by JSPS KAKENHI Grant Number JP19K03480.

\section{Preliminaries}
Let \((M,g)\)  be an \(n\)-dimensional Riemannian spin manifold with spinor bundle \(S_{1/2}\) and the complexified tangent bundle \(TM^{\mathbb{C}}\) . We consider the twisted Dirac operator on \(S_{1/2}\otimes TM^{\mathbb{C}}\)
\[D_{TM}=\sum_{k=1}^{n} (e_k\cdot \otimes \mathop{\mathrm{id}_{TM^\mathbb{C}}})\circ \nabla_{e_k}\]
where \(\nabla\) is the covariant derivative on \(S_{1/2}\otimes TM^{\mathbb{C}}\) and \(e_k\cdot\) is the Clifford multiplication by an orthonormal frame \(\{e_k\}\). 
With respect to  the \(\mathrm{Spin}(n)\) decomposition \(S_{1/2}\otimes TM^{\mathbb{C}} \cong S_{1/2} \oplus S_{3/2}\), we can write \(D_{TM}\) as the \(2\times 2\)-matrix
\[D_{TM}=
\begin{pmatrix}
\frac{2-n}{n}D&2P^\ast\\[1ex]
\frac{2}{n}P&R
\end{pmatrix}
\]
where \(P:\Gamma(S_{1/2})\to\Gamma(S_{3/2})\) is the twistor operator and \(P^\ast\) is the formal adjoint operator of \(P\). The operator \(R:\Gamma(S_{3/2})\to\Gamma(S_{3/2})\) is called the Rarita-Schwinger operator, which is a first order elliptic differential operator with conformal covariance. Furthermore, there is the $L^2$-decomposition \(\Gamma(S_{3/2})=\mathop{\mathrm{Ker}}P^\ast \oplus \mathop{\mathrm{Im}}P\) on a compact Riemannian spin manifold since twistor operator \(P\) is an overdetermined elliptic operator \cite{Be}.

Some Weitzenb\"ock formulas in \cite{HS} and \cite{Wa} gives us the following formulas about the Rarita-Schwinger operator on Einstein manifolds.

\begin{thm}[\cite{HS}]\label{HS}
Let \((M,g)\) be an $n$-dimensional compact Einstein spin manifold, then
\begin{enumerate}
\item \(R^2=(\frac{n-2}{n})^2(\Delta_{3/2}+\frac{1}{8}\mathrm{scal})\)\quad on \(\mathop{\mathrm{Im}}P\),
\item \(R^2=\Delta_{3/2}+\frac{n-8}{8n}\mathrm{scal}\)\quad on \(\mathop{\mathrm{Ker}}P^\ast\).
\end{enumerate}
Here, \(\Delta_{k/2}\) is the standard Laplacian on the bundle \(S_{k/2}\) introduced in \cite{SW}. For example, \(\Delta_{1/2}=\nabla^{\ast}\nabla+\frac{\mathrm{scal}}{8}=D^2-\frac{\mathrm{scal}}{8}\). Furthermore, we have equations about commutation,
\[R\circ P=\frac{n-2}{n}P\circ D,\quad P^\ast\circ R=\frac{n-2}{n}D\circ P^\ast,\]
\[\Delta_{1/2}\circ D=D\circ\Delta_{1/2},\quad
\Delta_{1/2}\circ P^\ast=P^\ast\circ\Delta_{3/2},\quad
\Delta_{3/2}\circ P=P\circ\Delta_{1/2},\quad
\Delta_{3/2}\circ R=R\circ\Delta_{3/2}.\]
In particular, \(\mathop{\mathrm{Im}}P\) and \(\mathop{\mathrm{Ker}}P^\ast\) are invariant under the action of the Rarita-Schwinger operator \(R\).
\end{thm}

Note that, on a homogeneous vector bundle over a compact symmetric space \(E=G\times_{\rho}V\), the standard Laplacian \(\Delta_{E}\) on the sections of \(E\) coincides with the Casimir operator of \(G\) \cite{SW}.

\section{The spectra of the Rarita-Schwinger operator on symmetric spaces}
Let \(M=G/K\) be a compact symmetric space and \(E\) be a homogeneous vector bundle over \(M\). Then, by Frobenius reciprocity, we can decompose the space of \(L^2\)-sections \(L^2(E)\) into the Hilbert sum
\[L^2(E)=\bigoplus_{\lambda\in\widehat{G}}
\mathrm{Hom}_{K}(V_\lambda, E)\otimes V_\lambda\]
where \(\widehat{G}\) is the set of equivalence classes of the irreducible representations of \(G\) 
and \(V_\lambda\) is the irreducible representation space of highest weight \(\lambda\). It follows from the Freudenthal's formula that the restriction of the standard Laplacian $\Delta_E$ on \(E\) to the above \(V_\lambda\) satisfies
\[\Delta_{E}\vert_{V_{\lambda}}=\langle \lambda+2\delta_{G}, \lambda \rangle\]
where \(\delta_{G}\) is half the sum of the positive roots on \(\mathfrak{g} =\mathrm{Lie}(G)\) and this inner product on weights comes from the Killing form on \(\mathfrak{g}\). 
 
 From now on, we consider an irreducible compact symmetric space $(G/K,g)$ with spin structure  (cf. \cite{CG2}). Then the metric $g$ is Einstein, and the space of the twistor spinors \(\mathop{\mathrm{Ker}}P\) coincides with the space of the real Killing spinors. Since a symmetric space admitting a real Killing spinor is locally conformally flat (cf. \cite{BMS}), if $G/K$ has nontrivial \(\mathop{\mathrm{Ker}}P\), then it is a space form of positive curvature, that is, \(S^n/\Gamma\) where \(\Gamma\) is a discrete group of $\mathrm{O}(n+1)$. Furthermore, the manifold is \(S^n\) or \(\mathbb{R}P^n\) because the fundamental group of symmetric space \(G/K\) is included in the center of \(G\). Note that \(\mathbb{R}P^n\) has a spin structure only for $n\equiv 3\mod 4$. In \cite {Bar}, it was proved that \(S^n\) and \(\mathbb{R}P^{4m-1}\) have real Killing spinors. They constitute $G$-modules whose highest weights are easily understood. As a result we know $G$-module structure of \(\mathop{\mathrm{Im}}P\) isomorphic to $L^2(S_{1/2})\ominus \mathop{\mathrm{Ker}}P$, and hence, its orthogonal complement \(\mathop{\mathrm{Ker}}P^{\ast}\) in $L^2(S_{3/2})$. In fact we will give irreducible decompositions of them on the standard sphere $S^n$ in the next subsection. When an irreducible compact symmetric space $G/K$ is not $S^n$ or $\mathbb{R}P^n$, the kernel of the twistor operator \(\mathop{\mathrm{Ker}}P\) is zero. Then the space \(\mathop{\mathrm{Im}}P\) in $L^2(S_{3/2})$ is isomorphic to $L^2(S_{1/2})$ as a $G$-module. By Frobenius reciprocity, we decompose $L^2(S_{1/2})$ and $L^2(S_{3/2})$ if we know the branching rule for $G/K$. By subtracting \(\mathop{\mathrm{Im}}P\) from $L^2(S_{3/2})$, we can have an irreducible decomposition of $\mathop{\mathrm{Ker}}P^\ast$. As examples, we will decompose \(\mathop{\mathrm{Im}}P\) and $\mathop{\mathrm{Ker}}P^\ast$ on $\mathbb{C}P^n$ and $\mathbb{H}P^n$ later. Thus we have known $G$-module structure of \(\mathop{\mathrm{Im}}P\) and $\mathop{\mathrm{Ker}}P^\ast$ on irreducible compact symmetric spaces. 

We calculate the eigenvalues of the square of the Rarita-Schwinger operator $R^2$. It follows from Theorem \ref{HS} that,  for $\psi$ in an irreducible summand $V_{\lambda}$ in \(\mathop{\mathrm{Im}}P\), 
\begin{equation*}
R^2\psi=\left(\frac{n-2}{n}\right)^2\big(\Delta_{3/2}+\frac{1}{8}\mathrm{scal}\big)\psi=\left(\frac{n-2}{n}\right)^2\big(\langle \lambda+2\delta_G,\lambda\rangle+\frac{1}{8}\mathrm{scal}\big)\psi.
\end{equation*}
\begin{rem}
Since we have $R^2P\phi=\left(\frac{n-2}{n}\right)^2PD^2\phi$ for spinor filed $\phi$, the eigenvalues for $R^2$ on \(\mathop{\mathrm{Im}}P\) correspond to the eigenvalues for $D^2$ up to a constant multiple. 
\end{rem}
When we take $\psi$ in an irreducible summand $V_{\lambda}$ in \(\mathop{\mathrm{Ker}}P^\ast\), Theorem \ref{HS} again allows us to get the eigenvalue, 
\begin{equation*}
R^2\psi=\big(\Delta_{3/2}+\frac{n-8}{8n}\mathrm{scal}\big)\psi=\big(\langle \lambda+2\delta_G,\lambda\rangle+\frac{n-8}{8n}\mathrm{scal}\big)\psi. 
\end{equation*}
We shall give some examples in the next subsections. For simplicity, an irreducible representation space \(V_\lambda\) is sometimes denoted by its highest weight \(\lambda\). In addition, \(k_j\) denotes a string of \(k\) with length \(j\). For example, \(((1/2)_{m-1},\pm 1/2)=(\underbrace{1/2,\dots,1/2}_{m-1},\pm 1/2)\).
\begin{rem}
The kernel of the Rarita-Schwinger operator $\mathop{\mathrm{Ker}}R$ is the space of the Rarita-Schwinger fields, which are important fermion fields from a view point of physics and geometry (cf. \cite{HS}). By using Theorem \ref{HS}, we can easily see that $\mathop{\mathrm{Ker}}R$ for compact $G/K$ with dimension $n>8$ is zero. Furthermore, the irreducible compact symmetric spaces admitting Rarita-Schwinger fields are classified in \cite{HS}, where \(\mathrm{Gr}_2(\mathbb{C}^4)\) and \(\frac{\mathrm{SO}(6)}{\mathrm{SO}(2)\times\mathrm{SO}(4)}\) are separated, but they are the same as a symmetric space.
\end{rem}

\subsection{the sphere}
The \(n\)-dimensional standard sphere \(S^n\) is recognized as the symmetric space \(S^n=\mathrm{Spin}(n+1)/\mathrm{Spin}(n)\). We put \(m=\lfloor (n+1)/2\rfloor\). The highest weight of \(\mathrm{Spin}(n+1)\) is written as
\(\lambda=(\lambda_1,\dots,\lambda_m)\) in \(\mathbb{Z}^m\) or \((\mathbb{Z}+1/2)^m\) which satisfies the dominant condition
\begin{align*}
&\lambda_1\geq \lambda_2\geq\dots\geq\lambda_{m-1}\geq\lambda_{m}\geq0 &&\mbox{for } n=2m,\\
&\lambda_1\geq \lambda_2\geq\dots\geq\lambda_{m-1}\geq\lvert\lambda_{m}\rvert &&\mbox{for } n=2m-1.
\end{align*}

The spinor bundle  \(S_{1/2}\) over the sphere \(S^n\) is  realized as the homogeneous vector bundle \(\mathrm{Spin}(n+1)\times_{\mathrm{Spin}(n)}V_{\lambda}\) with the highest weight 
\begin{align*}
&((1/2)_{m-2},1/2)\oplus((1/2)_{m-2},-1/2)&&\mbox{for } n=2m,\\
&((1/2)_{m-2},1/2) &&\mbox{for } n=2m-1
\end{align*}
of \(\mathrm{Spin}(n)\). By Frobenius reciprocity and branching rule (cf.\cite{Z}),
\(L^2(S_{1/2})\) is decomposed into
the irreducible $\mathrm{Spin}(n+1)$-modules with the highest weight \(\lambda=(\lambda_1,\dots,\lambda_m)\) which satisfy the dominant condition and
\begin{align*}
&\lambda_1\geq 1/2\geq\lambda_2\geq1/2\geq\dots\geq1/2\geq\lambda_{m-1}\geq \pm1/2\geq -\lambda_{m} &&\mbox{for } n=2m,\\
&\lambda_1\geq 1/2\geq \lambda_2\geq 1/2\geq \dots\geq1/2\geq\lambda_{m-1}\geq1/2\geq \lvert\lambda_{m}\rvert &&\mbox{for } n=2m-1.
\end{align*}
Especially, if \(\lambda\) satisfies this conditions, then \(\mathop{\mathrm{dim}}\mathop{\mathrm{Hom}_{\mathrm{Spin(n)}}}(V_\lambda, S_{1/2})=1\).
The highest weights which satisfy the above conditions are
\begin{align}
\lambda=
\begin{cases}
(k+1/2, (1/2)_{m-2}, 1/2)&\mbox{for } n=2m,\\
(k+1/2, (1/2)_{m-2}, \pm 1/2) &\mbox{for } n=2m-1,
\end{cases}\label{L2S1s}
\end{align}
where \(k\) is non-negative integer. Thus, we have the following decomposition 
\[L^2(S_{1/2})=\bigoplus_{\lambda}V_{\lambda}\]
where \(\lambda\) runs over the highest weights in (1).
On these irreducible components, the Dirac operator  satisfies
\begin{align}
D^2\vert_{V_\lambda}
&=\Delta_{\frac{1}{2}}\vert_{V_\lambda}+\frac{1}{8}\mathrm{scal}\notag\\
&=\langle \lambda+2\delta_{\mathrm{Spin}(n+1)}, \lambda\rangle+\frac{n(n-1)}{2} \notag
\\
&=\left(k+\frac{n}{2}\right)^2 \label{DEs}
\end{align}
because \(2\delta_{\mathrm{Spin}(n+1)}=(n-1, n-3, \dots, n+1-2m)\).

\begin{rem}
The inner product in the above equations satisfies \(\langle \mathbf{e}_{i}, \mathbf{e}_{j}\rangle=\delta_{ij}\) using the standard basis \(\{\mathbf{e}_i=(0,\dots,0,\overset{i}{1},0,\dots,0)\}\) of the dual of the maximal abelian subalgebra \(\mathfrak{h}^\ast\).
\end{rem}

Next, since the highest weight of \(TM^\mathbb{C}\) with respect to \(\mathrm{Spin}(n)\) is \((1,0,\dots,0)\), we get the irreducible decomposition
\[S_{1/2}\otimes TM^{\mathbb{C}}=
\begin{cases}
(3/2,(1/2)_{m-3},\pm1/2)\oplus (1/2,(1/2)_{m-3},\mp1/2) &\mbox{for } n=2m,\\
(3/2,(1/2)_{m-3},1/2)\oplus (1/2,(1/2)_{m-3},1/2) &\mbox{for } n=2m-1.
\end{cases}
\]
Therefore, 
\[
S_{3/2}=
\begin{cases}
(3/2,(1/2)_{m-3},1/2)\oplus (3/2,(1/2)_{m-3},-1/2) &\mbox{for } n=2m,\\
(3/2,(1/2)_{m-3},1/2) &\mbox{for } n=2m-1.
\end{cases}
\]
By Frobenius reciprocity and branching rule,  \(L^2(S_{3/2})\) is decomposed into
the irreducible representations of the highest weight \(\lambda=(\lambda_1,\dots,\lambda_m)\) which satisfies the dominant condition and
\begin{align*}
&\lambda_1\geq 3/2\geq\lambda_2\geq1/2\geq\dots\geq1/2\geq\lambda_{m-1}\geq \pm1/2\geq -\lambda_{m} &&\mbox{for } n=2m,\\
&\lambda_1\geq 3/2\geq \lambda_2\geq 1/2\geq\dots\geq1/2\geq\lambda_{m-1}\geq1/2\geq \lvert\lambda_{m}\rvert &&\mbox{for } n=2m-1.
\end{align*}
Especially, if \(\lambda\) satisfies this conditions, then \(\mathop{\mathrm{dim}}\mathop{\mathrm{Hom}}_{\mathrm{Spin}(n)}(V_\lambda, S_{3/2})=1\).
The highest weights which satisfy the above conditions are
\[
\lambda=
\begin{cases}
(l+3/2, 3/2, (1/2)_{m-2}), (l+3/2, 1/2, (1/2)_{m-2})&\mbox{for } n=2m,\\
(l+3/2, 3/2, (1/2)_{m-3}, \pm 1/2),  (l+3/2, 1/2, (1/2)_{m-3}, \pm 1/2)&\mbox{for } n=2m-1,
\end{cases}
\]
where  \(l\) is non-negative integer. 

On a compact \(n\)-dimensional Einstein spin manifold, the space of the twistor spinors coincides with the eigenspace of Dirac operator with the eigenvalue \(c_{0}\) which satisfies
\[c_{0}^2=\frac{n}{4(n-1)}\mathrm{scal}.\quad (\mbox{cf. }\cite{Fr})\]

By the equation (\ref{DEs}) and this fact , we have
\[
\mathop{\mathrm{Ker}}P=
\begin{cases}
(1/2, 1/2, 1/2, \dots, 1/2)&\mbox{for } n=2m,\\
(1/2, 1/2, 1/2, \dots, 1/2, \pm 1/2)&\mbox{for } n=2m-1.
\end{cases}
\]

Therefore, considering the highest weights of \(L^2(S_{1/2})\) in (\ref{L2S1s}), we have
\begin{subequations}
\begin{align}
\mathop{\mathrm{Im}}P&=
\begin{cases}
\bigoplus_{k\geq 1}(k+1/2, 1/2, (1/2)_{m-3}, 1/2)&\mbox{for } n=2m,\\
\bigoplus_{k\geq 1}(k+1/2, 1/2, (1/2)_{m-3}, \pm 1/2)&\mbox{for } n=2m-1,
\end{cases}\label{Ims}
\\
\mathop{\mathrm{Ker}}P^\ast&=
\begin{cases}
\bigoplus_{l\geq0}(l+3/2, 3/2, (1/2)_{m-3}, 1/2)&\mbox{for } n=2m,\\
\bigoplus_{l\geq0}(l+3/2, 3/2, (1/2)_{m-3}, \pm 1/2)&\mbox{for } n=2m-1.
\end{cases}\label{Kers}
\end{align}
\end{subequations}
Here, \(k\) is positive integer and \(l\) is non-negative integer. 

From the eigenvalues of $D^2$, we get the eigenvalues of the square of the Rarita-Schwinger operator on the irreducible components in \(\mathop{\mathrm{Im}}P\)
\[\left(\frac{n-2}{n}\right)^2\left(k+\frac{n}{2}\right)^2\quad\mbox{on (\ref{Ims}).}\]
On the irreducible components \(V_\lambda\) of \(\mathop{\mathrm{Ker}}P^\ast\subset\Gamma(S_{3/2})\) in (\ref{Kers}),
\begin{align*}
R^2|_{V_{\lambda}}&=\Delta_{3/2}|_{V_{\lambda}}+\frac{n-8}{8n}\mathrm{scal}\\
&=\langle \lambda+2\delta_{\mathrm{Spin}(n+1)}, \lambda\rangle+\frac{(n-8)(n-1)}{8}\\
&=\left(l+\frac{n+2}{2}\right)^2.
\end{align*}
\begin{thm}[\cite{Br},\cite{BS}]\label{sphere}
The eigenvalues of the square of the Rarita-Schwinger operator on the sphere are
\begin{enumerate}
\item \(\mbox{on }\mathop{\mathrm{Im}}P\),
\begin{align*}
&\left(\frac{n-2}{n}\right)^2\left(k+\frac{n}{2}\right)^2 \mbox{ with multi. } 2^{\lfloor \frac{n+1}{2}\rfloor}
\begin{pmatrix}k+n-1\\ k \end{pmatrix} \quad(k=1,2,\dots),
\end{align*}
\item \(\mbox{on }\mathop{\mathrm{Ker}}P^\ast\),
\begin{align*}
&\left(l+1+\frac{n}{2}\right)^2
\mbox{ with multi. } 
2^{\lfloor\frac{n+1}{2}\rfloor}(n-2)\frac{(l+n+1)(l+1)}{(l+n)(l+2)}\binom{l+n}{l+1}
\quad (l=0,1,2,\dots).
\end{align*}
\end{enumerate}
\end{thm}

\subsection{the complex projective space}
In this subsection, we shall calculate the case of the complex projective spaces with Fubini-Study metric, \(\mathbb{C}P^n=\mathrm{SU}(n+1)/\mathrm{S}(\mathrm{U}(n)\times\mathrm{U}(1))\). If \(n\) is even, the complex projective space \(\mathbb{C}P^n\) has no spin structure, so we consider only the case that \(n\) is odd and \(n\geq3\).

Since the isotropy representation is an \((n+1)\)-fold covering homomorphism \(\mathrm{S}(\mathrm{U}(n)\times\mathrm{U}(1))\ni (A,a)\mapsto a^{-1}A\in \mathrm{U}(n)\), we get the correspondence between the highest weight of \(\mathrm{U}(n)\) and \(\mathrm{S}(\mathrm{U}(n)\times\mathrm{U}(1))\) for the irreducible \(\mathrm{U}(n)\)-modules,
\begin{equation}\label{transhw}
(\mu_1,\dots,\mu_n)_{\mathrm{U}(n)}
=\left(\mu_1+\sum_{i=1}^n\mu_i\dots,\mu_n+\sum_{i=1}^n\mu_i\right)_{\mathrm{S}(\mathrm{U}(n)\times\mathrm{U}(1))}
\end{equation}
where we use a notation in \cite{CG} to write weights for \(\mathrm{S}(\mathrm{U}(n)\times\mathrm{U}(1))\). For example, the highest weight of the canonical bundle $\Lambda^{n,0}$ is (\((-1)_n)_{\mathrm{U}(n)}=((-(n+1))_{n})_{\mathrm{S}(\mathrm{U}(n)\times\mathrm{U}(1))}\). Through the correspondence, we get the information on the highest weight for \(\mathrm{S}(\mathrm{U}(n)\times\mathrm{U}(1))\). Indeed we know the dominant condition for the \(\mathrm{S}(\mathrm{U}(n)\times\mathrm{U}(1))\)-module \(\mu'=(\mu'_1,\dots,\mu'_{n})\) is
\[\mu'_1\geq\cdots\geq\mu'_{n}.\] 
We also know there exists the square root of the canonical bundle \(\sqrt{\Lambda^{n,0}}\) on $\mathbb{C}P^n$, which gives a spin structure, because its highest weight \(\left(\Bigl(-\frac{n+1}{2}\Bigr)_{n}\right)_{\mathrm{S}(\mathrm{U}(n)\times\mathrm{U}(1))}=\left(\Bigl(-\frac{1}{2}\Bigr)_{n}\right)_{\mathrm{U}(n)}\) is an integral weight for \(\mathrm{S}(\mathrm{U}(n)\times\mathrm{U}(1))\) when \(n\) is odd.
 
The spinor bundle is decomposed by the action of \(\mathrm{S}(\mathrm{U}(n)\times\mathrm{U}(1))\) into
\[S_{1/2}=\bigoplus_{k=0}^n S_{1/2}(k)\]
where \(S_{1/2}(k)\) is the irreducible bundle of the highest weight
\[\left(\Bigl(\frac{n+1}{2}-k\Bigr)_{n-k},
\Bigl(\frac{n-1}{2}-k\Bigr)_{k}\right).\quad(\mbox{cf. }\cite{SS})\]

Using the same method of the case of the sphere, we get the eigenvalues of the Rarita-Schwinger operator on the complex projective space. By Frobenius reciprocity and branching rule in \cite{CG}, \cite{SS},  \(L^2(S_{1/2}(k))\) is decomposed into the direct sum of
the irreducible \(\mathrm{SU}(n+1)\)-modules of the highest weight $\lambda=(\lambda_1,\dots,\lambda_n)$ in $\mathbb{Z}^n$ which satisfies the dominant condition for \(\mathrm{SU}(n+1)\)
\[\lambda_1\geq\lambda_2\geq\cdots\geq\lambda_n\geq0\]
and the condition that there exists 
$(\nu_1,\dots, \nu_n)$ in $\mathbb{Z}^n$ such that

\[l:=\sum_{i=1}^n(\lambda_i-\nu_i)\]
\[\frac{n+1}{2}-k=\nu_i-l \quad (i=1,\dots,n-k),\qquad
\frac{n-1}{2}-k=\nu_i-l \quad (i=n-k+1,\dots,n)\]
\[\lambda_1\geq\nu_1\geq\lambda_2\geq\nu_2\geq\dots\geq\lambda_n\geq\nu_n\geq 0.\]
Especially, if \(\lambda\) satisfies this conditions, then \(\mathop{\mathrm{dim}}\mathop{\mathrm{Hom}_{\mathrm{S}(\mathrm{U}(n)\times\mathrm{U}(1))}}(V_\lambda, S_{1/2}(k))=1\). Thus the highest weight \(\lambda\) in \(L^2(S_{1/2}(k))\) is one of the followings.
\begin{subequations}If \(k=1,\dots,n-1\), then
\begin{equation}
\lambda(k,\epsilon,l)=\left(\frac{n+1}{2}-k+2l-\epsilon,\Bigl(\frac{n+1}{2}-k+l\Bigr)_{n-k-1},
\frac{n-1}{2}-k+l+\epsilon,\Bigl(\frac{n-1}{2}-k+l\Bigr)_{k-1}\right) \label{L2S1cp1}
\end{equation}
where \(\epsilon\in\{0,1\}\) and \(l\geq \max\{\epsilon,-\frac{n-1}{2}+k\}\).
If \(k=0\), then
\begin{equation}
\lambda(0,0,l)=\left(\frac{n+1}{2}+2l,\Bigl(\frac{n+1}{2}+l\Bigr)_{n-1}\right) \label{L2S1cp2}
\end{equation}
where \(l\geq 0\). 
If \(k=n\), then
\begin{equation}
\lambda(n,1,l)=\left(\frac{n-1}{2}-n+2l,\Bigl(\frac{n-1}{2}-n+l\Bigr)_{n-1}\right) \label{L2S1cp3}
\end{equation}
where \(l\geq \frac{n+1}{2}\).
\end{subequations}

On these irreducible components, the Dirac operator  satisfies
\begin{align}\label{DEcp}
D^2\vert_{V_\lambda}
&=\Delta_{\frac{1}{2}}\vert_{V_\lambda}+\frac{1}{8}\mathrm{scal}\notag\\
&=\langle \lambda+2\delta_{\mathrm{SU}(n+1)}, \lambda\rangle+\frac{n(n+1)}{4}\notag\\
&=
\begin{cases}
 (l+n-k)(2l+n+1-2\epsilon)   &\mbox{on (\ref{L2S1cp1})}\\
 (l+n)(2l+n+1)  &\mbox{on (\ref{L2S1cp2})} \\
 l(2l+n-1)  &\mbox{on (\ref{L2S1cp3})} 
\end{cases}
\end{align}
because \(2\delta_{\mathrm{SU}(n+1)}=(2n, 2(n-1), \dots, 2)\).

\begin{rem}
The inner product of the above equations satisfies \(\langle \mathbf{e}_{i}, \mathbf{e}_{j}\rangle=\frac{(n+1)\delta_{ij}-1}{n+1}\) using the basis \(\{\mathbf{e}_i=(0,\dots,0,\overset{i}{1},0,\dots,0)\}\) of the dual of the maximal abelian subalgebra \(\mathfrak{h}^\ast\).
\end{rem}

Next, we consider the decomposition of  \(S_{1/2}(k)\otimes TM^\mathbb{C}\). The tangent bundle $TM^{\mathbb{C}}$ is the direct sum of $(1,0)$-part \(T^{1,0}M\) and $(0,1)$-part  \(T^{0,1}M\). Since the highest weight of \(T^{1,0}M\) and \(T^{0,1}M\) is \((1,0,\cdots0)\) and \((0,\cdots,0,-1)\) with respect to the action of \(\mathrm{U}(n)\), we have the decomposition with respect to \(\mathrm{U}(n)\) and also \(\mathrm{S}(\mathrm{U}(n)\times\mathrm{U}(1))\) by the correspondence (\ref{transhw}),
\begin{align*}
S_{1/2}(k)\otimes TM^\mathbb{C}=
\begin{cases}
S_{3/2}^+(k)\oplus S_{1/2}(k-1)\oplus S_{1/2}(k+1)\oplus S_{3/2}^-(k)
 &\mbox{for } k=1,\dots,n-1\\
S_{3/2}^+(0)\oplus S_{1/2}(1)&\mbox{for } k=0\\
S_{1/2}(n-1)\oplus S_{3/2}^-(n)&\mbox{for } k=n\\
\end{cases}
\end{align*}
Here, we put 
\[S_{3/2}^+(k)=
\left(\frac{n+5}{2}-k,\Bigl(\frac{n+3}{2}-k\Bigr)_{n-k-1}, \Bigl(\frac{n+1}{2}-k\Bigr)_{k}\right),\]
\[S_{3/2}^-(k)=
\left(\Bigl(\frac{n-1}{2}-k\Bigr)_{n-k}, \Bigl(\frac{n-3}{2}-k\Bigr)_{k-1}, \frac{n-5}{2}-k\right).\]
Thus, we get
\[S_{3/2}=\bigoplus_{k=1}^{n-1}S_{1/2}(k)\bigoplus_{k=0}^{n-1}S_{3/2}^+(k)\bigoplus_{k=1}^{n}S_{3/2}^-(k).\]

As the case of the spinor bundle , \(L^2(S_{3/2}^+(k))\) is decomposed into the irreducible \(\mathrm{SU}(n+1)\)-modules with the highest weights $\lambda=(\lambda_1,\dots,\lambda_n)$ in $\mathbb{Z}^n$ which satisfies the condition that $\lambda$ is dominant and there exists $(\nu_1,\dots, \nu_n)$ in $\mathbb{Z}^n$ such that
\begin{align*}
l&=\sum_{i=1}^n(\lambda_i-\nu_i),\\
\frac{n+5}{2}-k&=\nu_1-l ,\\
\frac{n+3}{2}-k&=\nu_i-l \quad (i=2,\dots,n-k),\\
\frac{n+1}{2}-k&=\nu_i-l \quad (i=n-k+1,\dots,n), \\
\lambda_1\geq\nu_1\geq&\lambda_2\geq\nu_2\geq\dots\geq\lambda_n\geq\nu_n\geq 0.
\end{align*}
Especially, if \(\lambda\) satisfies this conditions, then \(\mathop{\mathrm{dim}}\mathop{\mathrm{Hom}_{\mathrm{S}(\mathrm{U}(n)\times\mathrm{U}(1))}}(V_\lambda, S^+_{3/2}(k))=1\). Then, \(\lambda\) is one of the followings. If $k=1,\dots, n-2$, then
\begin{align*}
\lambda^+(k,\epsilon,\epsilon_1,l)
&=\left(\frac{n+5}{2}-k+2l-\epsilon-\epsilon_1,\frac{n+3}{2}-k+l+\epsilon_1,\Bigl(\frac{n+3}{2}-k+l\Bigr)_{n-k-2},\right.\\
&\hspace{17em}\left.\frac{n+1}{2}-k+l+\epsilon,\Bigl(\frac{n+1}{2}-k+l\Bigr)_{k-1}\right)
\end{align*}
where \(\epsilon,\epsilon_1\in\{0,1\}\) and 
\(l\geq \max\{\epsilon+\epsilon_1,-\frac{n+1}{2}+k\}\).
If \(k=0\), then
\[\lambda^+(0,0,\epsilon_1,l)=\left(\frac{n+5}{2}+2l-\epsilon_1,\frac{n+3}{2}+l+\epsilon_1,
\Bigl(\frac{n+3}{2}+l\Bigr)_{n-2}\right).\]
where \(\epsilon_1\in\{0,1\}\) and \(l\geq\epsilon_1\). 
If \(k=n-1\), then
\[\lambda^+(n-1,1,\epsilon_1,l)=\left(\frac{n+5}{2}-(n-1)+2l-(1+\epsilon_1),\frac{n+1}{2}-(n-1)+l+(1+\epsilon_1),\Bigl(\frac{n+1}{2}-(n-1)+l\Bigr)_{n-2}\right).\]
where \(\epsilon_1\in\{-1,0,1\}\) and \(l\geq \max\{1+\epsilon_1,\frac{n-3}{2}\}\). 
\begin{rem}\label{rem+}
The $\mathrm{SU}(n+1)$-module $\lambda^+(n-1,1,-1,l-2)$ is isomorphic to $\lambda(0,0,l)$ in \eqref{L2S1cp2}, hence $\bigoplus_l\lambda^+(n-1,1,-1,l-2)\cong L^2(S_{1/2}(0))$. 
\end{rem}

Next, \(L^2(S_{3/2}^-(k))\) is decomposed into the irreducible representations with the highest weights $\lambda=(\lambda_1,\dots,\lambda_n)$ in $\mathbb{Z}^n$ which satisfies the condition that $\lambda$ is dominant there exists $(\nu_1,\dots, \nu_n)$ in $\mathbb{Z}^n$ such that
\begin{align*}
l&=\sum_{i=1}^n(\lambda_i-\nu_i),\\
\frac{n-1}{2}-k&=\nu_i-l \quad (i=1,\dots,n-k),\\
\frac{n-3}{2}-k&=\nu_i-l \quad (i=n-k+1,\dots,n-1), \\
\frac{n-5}{2}-k&=\nu_n-l , \\
\lambda_1\geq\nu_1\geq&\lambda_2\geq\nu_2\geq\dots\geq\lambda_n\geq\nu_n\geq 0.
\end{align*}
Especially, if \(\lambda\) satisfies this conditions, then \(\mathop{\mathrm{dim}}\mathop{\mathrm{Hom}_{\mathrm{S}(\mathrm{U}(n)\times\mathrm{U}(1))}}(V_\lambda, S^-_{3/2}(k))=1\). Then, \(\lambda\) is one of the followings. If $k=2,\dots,n-1$, then 
\begin{align*}
\lambda^-(k,\epsilon,\epsilon_2,l)
=&\left(\frac{n-1}{2}-k+2l-\epsilon-\epsilon_2,\Bigl(\frac{n-1}{2}-k+l\Bigr)_{n-k-1},\frac{n-3}{2}-k+l+\epsilon,\right.\\
&\left.\hspace{15em}\Bigl(\frac{n-3}{2}-k+l\Bigr)_{k-2},\frac{n-5}{2}-k+l+\epsilon_2\right)
\end{align*}
where \(\epsilon,\epsilon_2\in\{0,1\}\) and 
\(l\geq \max\{\epsilon+\epsilon_2,-\frac{n-5}{2}+k\}\).
If \(k=1\), then
\[\lambda^-(1,0,\epsilon_2,l)=\left(\frac{n-1}{2}-1+2l-\epsilon_2,\Bigl(\frac{n-1}{2}-1+l\Bigr)_{n-2},
\frac{n-5}{2}-1+l+\epsilon_2\right)\]
where \(\epsilon_2\in\{0,1,2\}\) and \(l\geq \max\{\epsilon_2,\frac{7-n}{2}\}\). If \(k=n\), then
\[\lambda^-(n,1,\epsilon_2,l)=\left(\frac{n-3}{2}-n+2l-\epsilon_2,\Bigl(\frac{n-3}{2}-n+l\Bigr)_{n-2},\frac{n-5}{2}-n+l+\epsilon_2 \right)\]
where \(\epsilon_2\in\{0,1\}\) and \(l\geq\frac{n+5}{2}\). 

\begin{rem}\label{rem-}
We know $\bigoplus_l \lambda^-(1,0,2,l)\cong L^2(S_{1/2}(n))$ as a $\mathrm{SU}(n+1)$-module. 
\end{rem}

We shall give decomposition of $\mathop{\mathrm{Im}} P$ and $\mathop{\mathrm{Ker}} P^{\ast}$. Because of $\mathop{\mathrm{Ker}} P=\{0\}$, we have
\begin{align*}
\mathop{\mathrm{Im}}P\cong&L^2(S_{1/2})=
\bigoplus_{k=0}^{n}\bigoplus_{\epsilon,l}\lambda(k,\epsilon,l)
\end{align*}
where $l\ge \max \{\epsilon, -\frac{n-1}{2}+k\}$ and $\epsilon$ runs as follows (i) for  \(k=0\), \(\epsilon=0\)  (ii)  for \(k=1,\dots,n-1\), \(\epsilon=0,1\)  (iii)  for \(k=n\), \(\epsilon=1\). According to Remark \ref{rem+} and \ref{rem-}, we can decompose $\mathop{\mathrm{Ker}} P^{\ast}= L^2(S_{3/2})\ominus \mathop{\mathrm{Im}}P$. We divide it into two parts, 
\begin{equation*}
\mathop{\mathrm{Ker}}P^\ast=(\mathop{\mathrm{Ker}}P^\ast)^+\oplus(\mathop{\mathrm{Ker}}P^\ast)^-.
\end{equation*}
The first part is 
\begin{equation*}
(\mathop{\mathrm{Ker}}P^\ast)^+\cong
\bigoplus_{k=0}^{n-1}\bigoplus_{\epsilon,\epsilon_{1},l}\lambda^+(k,\epsilon,\epsilon_1,l),
\end{equation*}
where $l\geq \max\{\epsilon+\epsilon_1,-\frac{n+1}{2}+k\}$, and $\epsilon$, $\epsilon_1$ run as follows: (i) for $k=0$, $\epsilon=0$, $\epsilon_1=0,1$, (ii) for $k=1,\dots,n-2$, $\epsilon=0,1$, $\epsilon_1=0,1$, (iii) for $k=n-1$, $\epsilon=1$, $\epsilon_1=0,1$. The second part is 
\begin{equation*}
(\mathop{\mathrm{Ker}}P^\ast)^-\cong
\bigoplus_{k=1}^{n}\bigoplus_{\epsilon,\epsilon_{2},l}\lambda^+(k,\epsilon,\epsilon_2,l),
\end{equation*}
where $l\geq \max\{\epsilon+\epsilon_2,-\frac{n-5}{2}+k\}$, and $\epsilon$, $\epsilon_2$ run as follows: (i) for $k=1$, $\epsilon=0$, $\epsilon_2=0,1$, (ii) for $k=2,\dots,n-2$, $\epsilon=0,1$, $\epsilon_2=0,1$.  (iii) for $k=n$, $\epsilon=1$, $\epsilon_2=0,1$.

Finally, we obtain the following eigenvalues by the same calculus as the case of the sphere.

\begin{thm}\label{cplx}
The eigenvalues of the square of the Rarita-Schwinger operator on the odd dimensional complex projective space $\mathbb{C}P^n$ are
\begin{enumerate}
\item \(\mbox{on }\lambda(k,\epsilon,l)\subset\mathop{\mathrm{Im}}P\),
\begin{align*}
&\left(\frac{n-1}{n}\right)^2  (l+n-k)(2l+n+1-2\epsilon),
\end{align*}
\item \(\mbox{on }\lambda^+(k,\epsilon,\epsilon_{1},l)\subset\mathop{\mathrm{Ker}}P^\ast\),
\begin{align*}
&((l+1)+n-k-\epsilon_1)(2(l+1)+n+1-2\epsilon)-(n+1)(1-\epsilon_1),
\end{align*}
\item \(\mbox{on }\lambda^-(k,\epsilon,\epsilon_{2},l)\subset\mathop{\mathrm{Ker}}P^\ast\),
\begin{align*}
&((l-1)+n-k+(1-\epsilon_2))(2(l-1)+n+1-2\epsilon)-(n+1)\epsilon_2.
\end{align*}
\end{enumerate}
\end{thm}

\begin{rem}
The dissertation of U.Semmelmann \cite{S} gave the eigenvalues of the square of the twisted Dirac operator \(D_{TM}\) on the odd dimensional complex projective spaces.
\end{rem}

\subsection{the quaternionic projective space}
Next, we shall consider the case of the quaternionic projective spaces.
The quaternionic projective space \(\mathbb{H}P^n\) is the symmetric space \(\mathrm{Sp}(n+1)/(\mathrm{Sp}(1)\times \mathrm{Sp}(n))\) with a spin structure. We calculate the eigenvalues for \(n\geq2\) as the previous subsections. First, the spinor bundle is decomposed into 
\[S_{1/2}=\bigoplus_{k=0}^n S_{1/2}(k)\]
where \(S_{1/2}(k)\) is the irreducible bundle of the highest weight
\[(k,1_{n-k},
0_{k})\]
with respect to the action of \(\mathrm{Sp}(1)\times \mathrm{Sp}(n)\). Here, the first component and the other components are parametrized by the highest weight for \(\mathrm{Sp}(1)\) and for \(\mathrm{Sp}(n)\), respectively.

By branching rule in \cite{T}, \(L^2(S_{1/2}(k))\) is decomposed into the direct sum of the irreducible \(\mathrm{Sp}(n+1)\)-modules of the highest weight \(\lambda=(\lambda_1,\dots,\lambda_{n+1})\) in $\mathbb{Z}^{n+1}$ which satisfies the dominant condition for \(\mathrm{Sp}(n+1)\) 
\[\lambda_1\geq\lambda_2\geq\cdots\geq \lambda_{n+1}\geq0\]
and for \(\mu=(\mu_1,\dots,\mu_{n+1})=(k, 1_{n-k}, 0_{k})\),
\[\lambda_{i}\geq\mu_{i+1}\geq\lambda_{i+2}\ (1\leq i\leq n-1),\quad \lambda_{n}\geq\mu_{n+1}.\]
In addition, the multiplicity \(m(\lambda)\) of \(\lambda\)  is the coefficient of \(X^{\mu_1+1}\) in
\[(X-X^{-1})^{-n}\prod_{i=1}^{n+1}(X^{l_i+1}-X^{-(l_i+1)})\]
where 
\begin{align*}
&l_1=\lambda_1-\max\{\lambda_2,\mu_2\},\\
&l_i=\min\{\lambda_i,\mu_i\}-\max\{\lambda_{i+1},\mu_{i+1}\}\ (1\leq i\leq n),\\
&l_{n+1}=\min\{\lambda_{n+1},\mu_{n+1}\}.
\end{align*}

For simplicity, we put \(\lambda(m,l,\epsilon)=(m+l,m,1_{n-1-l-\epsilon},0_{l+\epsilon})\). 
The highest weight \(\lambda\) satisfying the above condition is one of the followings.
\begin{align*}
\lambda=
\begin{cases}
\lambda(m,k,1) & \mbox{ for } 0\leq k\leq n-2\\
\lambda(m,k+1,-1) & \mbox{ for } 0\leq k\leq n-1\\
\lambda(m,k-1,1) & \mbox{ for } 1\leq k\leq n-1\\
\lambda(m,k,-1) & \mbox{ for } 1\leq k\leq n\\
\lambda(0,n,-1) & \mbox{ for } k=n-1,n
\end{cases}
\end{align*}
where \(m\) is positive integer and all of the multiplicities are 1. Therefore, we obtain the decomposition
\[L^2(S_{1/2})=
\bigoplus_{l=1}^{n}\bigoplus_{m\geq1}2\lambda(m,l,-1)
\oplus\bigoplus_{l=0}^{n-2}\bigoplus_{m\geq1}2\lambda(m,l,1)
\oplus 2\lambda(0,n,-1).\]

On these irreducible components \(\lambda(m,l,\epsilon)\), the Dirac operator  satisfies
\begin{align*}
D^2\vert_{V_\lambda}
&=\Delta_{\frac{1}{2}}\vert_{V_\lambda}+\frac{1}{8}\mathrm{scal}\\
&=\langle \lambda+2\delta_{\mathrm{Sp}(n+1)}, \lambda\rangle+\frac{n}{4}\\
&=\frac{1}{2(n+2)}(n+m-\epsilon)(n+m+l+1+\epsilon)
\end{align*}
because \(2\delta_{\mathrm{Sp}(n+1)}=(2(n+1), 2n, \dots, 2)\).
\begin{rem}
The inner product of the above equations satisfies \(\langle \mathbf{e}_{i}, \mathbf{e}_{j}\rangle=\frac{1}{4(n+2)}\delta_{ij}\) using the basis \(\{\mathbf{e}_i=(0,\dots,0,\overset{i}{1},0,\dots,0)\}\) of the dual of the maximal abelian subalgebra \(\mathfrak{h}^\ast\).
\end{rem}

Next, considering that the highest weight of \(\mathrm{Sp}(1)\times\mathrm{Sp}(n)\)-representation \(TM^\mathbb{C}\) is \((1,1,0,\dots,0)\), we have the irreducible decomposition
\begin{align*}
S_{1/2}\otimes TM^{\mathbb{C}}=&
\bigoplus_{k=0}^{n-2} (k,2,1_{n-k-2},0_{k+1})
\oplus\bigoplus_{k=1}^{n} (k,2,1_{n-k},0_{k-1})\\
&\oplus\bigoplus_{k=0}^{n-2} (k,1_{n-k-2},0_{k+2})
\oplus\bigoplus_{k=2}^{n+1} (k,1_{n-k+2},0_{k-2})
\oplus\bigoplus_{k=0}^{n-1} (k,1_{n-k},0_{k})
\oplus\bigoplus_{k=1}^{n} (k,1_{n-k},0_{k}).
\end{align*}
Therefore, 
\begin{align*}
S_{3/2}=&
\bigoplus_{k=0}^{n-2} (k,2,1_{n-k-2},0_{k+1})
\oplus\bigoplus_{k=1}^{n} (k,2,1_{n-k},0_{k-1})\\
&\oplus\bigoplus_{k=0}^{n-2} (k,1_{n-k-2},0_{k+2})
\oplus\bigoplus_{k=2}^{n+1} (k,1_{n-k+2},0_{k-2})
\oplus\bigoplus_{k=1}^{n-1} (k,1_{n-k},0_{k}).
\end{align*}
We denote by $L^2(\nu)$ the space of the $L^2$-sections of the vector bundle with the highest weight $\nu$ in the above decomposition. By branching rule, the space $L^2(\nu)$ is decomposed into the direct sum of irreducible $\mathrm{Sp}(n+1)$-modules. 

For simplicity, we put \(\tilde{\lambda}(m,l,\epsilon)=(m+l,m,2,1_{n-2-l-\epsilon},0_{l+\epsilon})\). The highest weight \(\lambda=(\lambda_1,\dots,\lambda_{n+1})\) of the irreducible summands in $L^2((k,2,1_{n-k-2},0_{k+1}))$ is
\begin{align*}
\lambda=
\begin{cases}
\tilde{\lambda}(m',k+1,0) & \mbox{ for } 0\leq k\leq n-3\\
\tilde{\lambda}(m',k-1,2) & \mbox{ for } 1\leq k\leq n-3\\
\tilde{\lambda}(m',k,0) & \mbox{ for } 0\leq k\leq n-2\\
\tilde{\lambda}(m',k,2) & \mbox{ for } 0\leq k\leq n-4\\
\lambda(m,k+2,-1) & \mbox{ for } 0\leq k\leq n-2\\
\lambda(m',k,1) & \mbox{ for } 0\leq k\leq n-3\\
\lambda(m,k,1) & \mbox{ for } 1\leq k\leq n-2\\
\lambda(m',k-2,3) & \mbox{ for } 2\leq k\leq n-2\\
\lambda(m,k+1,-1) & \mbox{ for } 0\leq k\leq n-2\\
\lambda(m,k+1,1) & \mbox{ for } 0\leq k\leq n-3\\
\lambda(m',k-1,1) & \mbox{ for } 1\leq k\leq n-2\\
\lambda(m',k-1,3) & \mbox{ for } 1\leq k\leq n-3\\
\lambda(0,n,-1) & \mbox{ for } k=n-2
\end{cases}
\end{align*}
where \(m\geq1\), \(m'\geq2\) and all of the multiplicities are $1$.

For $L^2((k,2,1_{n-k},0_{k-1}))$, the highest weight $\lambda$ of the irreducible summands is 
\begin{align*}
\lambda=
\begin{cases}
\tilde{\lambda}(m',k+1,-2) & \mbox{ for } 1\leq k\leq n-1\\
\tilde{\lambda}(m',k-1,0) & \mbox{ for } 1\leq k\leq n-1\\
\tilde{\lambda}(m',k,-2) & \mbox{ for } 2\leq k\leq n\\
\tilde{\lambda}(m',k,0) & \mbox{ for } 1\leq k\leq n-2\\
\lambda(m,k+2,-3) & \mbox{ for } 1\leq k\leq n\\
\lambda(m',k,-1) & \mbox{ for } 1\leq k\leq n-1\\
\lambda(m,k,-1) & \mbox{ for } 1\leq k\leq n\\
\lambda(m',k-2,1) & \mbox{ for } 2\leq k\leq n\\
\lambda(m,k+1,-3) & \mbox{ for } 2\leq k\leq n\\
\lambda(m,k+1,-1) & \mbox{ for } 1\leq k\leq n-1\\
\lambda(m',k-1,-1) & \mbox{ for } 2\leq k\leq n\\
\lambda(m',k-1,1) & \mbox{ for } 1\leq k\leq n-1\\
\lambda(0,n+2,-3)  & \mbox{ for } k=n-1
\end{cases}
\end{align*}
where \(m\geq1\), \(m'\geq2\) and all of the multiplicities are $1$.

For $L^2((k,1_{n-k-2},0_{k+2}))$, 
\begin{align*}
\lambda=
\begin{cases}
\lambda(m,k,1) & \mbox{ for } 0\leq k\leq n-2\\
\lambda(m,k+1,1) & \mbox{ for } 0\leq k\leq n-3\\
\lambda(m,k-1,3) & \mbox{ for } 1\leq k\leq n-3\\
\lambda(m,k,3) & \mbox{ for } 0\leq k\leq n-4\\
\lambda(0,n-2,1) & \mbox{ for } k=n-3,n-2
\end{cases}
\end{align*}
where \(m\) is positive integer and all of the multiplicities are $1$.

For $L^2((k,1_{n-k+2},0_{k-2}))$, 
\begin{align*}
\lambda=
\begin{cases}
\lambda(m,k,-1) & \mbox{ for } 2\leq k\leq n\\
\lambda(m,k+1,-3) & \mbox{ for } 2\leq k\leq n+1\\
\lambda(m,k-1,-1) & \mbox{ for } 2\leq k\leq n+1\\
\lambda(m,k,-3) & \mbox{ for } 3\leq k\leq n+1\\
\lambda(0,n+2,-3) & \mbox{ for } k=n+1
\end{cases}
\end{align*}
where \(m\) is positive integer and all of the multiplicities are $1$.

Because of $\mathop{\mathrm{Ker}} P=\{0\}$, we obtain
\begin{align*}
\mathop{\mathrm{Im}}P&\cong L^2(S_{1/2})=
\bigoplus_{l=1}^{n}\bigoplus_{m\geq1}2\lambda(m,l,-1)
\oplus\bigoplus_{l=0}^{n-2}\bigoplus_{m\geq1}2\lambda(m,l,1)
\oplus 2\lambda(0,n,-1),
\\
\mathop{\mathrm{Ker}}P^\ast&=
\bigoplus_{l=1}^{n-2}\bigoplus_{m\geq2}2\tilde{\lambda}(m,l,0)
\oplus\bigoplus_{l=0}^{n-2}\bigoplus_{m\geq2}2\tilde{\lambda}(m,l,0)
\oplus\bigoplus_{l=0}^{n-4}\bigoplus_{m\geq2}2\tilde{\lambda}(m,l,2)
\oplus\bigoplus_{l=2}^{n}\bigoplus_{m\geq2}2\tilde{\lambda}(m,l,-2)\\
&\qquad\oplus\bigoplus_{l=0}^{n-4}\bigoplus_{m\geq2}2\lambda(m,l,3)
\oplus\bigoplus_{l=0}^{n-4}\bigoplus_{m\geq1}2\lambda(m,l,3)
\oplus\bigoplus_{l=3}^{n+1}\bigoplus_{m\geq1}2\lambda(m,l,-3)
\oplus\bigoplus_{l=0}^{n+2}\bigoplus_{m\geq1}2\lambda(m,l,-3)\\
&\qquad\oplus\bigoplus_{l=0}^{n-3}\bigoplus_{m\geq2}2\lambda(m,l,1)
\oplus\bigoplus_{l=0}^{n-2}\bigoplus_{m\geq2}2\lambda(m,l,1)
\oplus\bigoplus_{l=1}^{n-2}\bigoplus_{m\geq2}4\lambda(m,l,1)\\
&\qquad\oplus\bigoplus_{l=1}^{n-1}\bigoplus_{m\geq2}2\lambda(m,l,-1)
\oplus\bigoplus_{l=2}^{n-1}\bigoplus_{m\geq1}2\lambda(m,l,-1)
\oplus\bigoplus_{l=1}^{n}\bigoplus_{m\geq1}2\lambda(m,l,-1)
\oplus\bigoplus_{l=2}^{n}\bigoplus_{m\geq1}2\lambda(m,l,-1)\\
&\qquad\oplus 2\lambda(0,n+2,-3)\oplus 2\lambda(0,n-2,1).
\end{align*}
Here, \(k\) is positive integer and \(l\) is non-negative integer.

We calculate the eigenvalue on each component. Then we have
\begin{thm}\label{qtn}
The eigenvalues of the square of the Rarita-Schwinger operator on the quaternionic projective space $\mathbb{H}P^n$ are
\begin{enumerate}
\item \(\mbox{on }\lambda(m,l,\epsilon)\subset\mathop{\mathrm{Im}}P\),
\begin{align*}
\left(\frac{2n-1}{2n}\right)^2 \frac{1}{2(n+2)}(n+m-\epsilon)(n+m+l+1+\epsilon),
\end{align*}
\item \(\mbox{on }\tilde{\lambda}(m,l,\epsilon)\subset\mathop{\mathrm{Ker}}P^\ast\), 
\begin{align*}
\frac{1}{2(n+2)}\bigl((n+m-\epsilon)(n+m+l+1+\epsilon)+\frac{1}{2}(\epsilon-2)(\epsilon+2)\bigr),
\end{align*}
\item \(\mbox{on }\lambda(m,l,\epsilon)\subset\mathop{\mathrm{Ker}}P^\ast\),
\begin{align*}
\frac{1}{2(n+2)}\bigl((n+m-\epsilon)(n+m+l+1+\epsilon)-(n+1/2) +\frac{1}{2}(\epsilon-2)(\epsilon+2)\bigr).
\end{align*}
\end{enumerate}
\end{thm}

\vspace{1\baselineskip}

Yasushi Homma, 
Department of Mathematics, Faculty of science and engineering, Waseda University, 3-4-1 Ohkubo, Shinjuku-ku, Tokyo, 169-8555, JAPAN.\\
E-mail address: homma\_yasushi@waseda.jp\\

Takuma Tomihisa, 
Department of Pure and applied Mathematics, Graduate school of fundamental science and engineering Waseda University, 3-4-1 Ohkubo, Shinjuku-ku, Tokyo, 169-8555, JAPAN.\\
E-mail address: taku-tomihisa@akane.waseda.jp

\end{document}